\documentclass[draft,reqno]{amsart}

\usepackage{amsmath,amsthm,amssymb,amsfonts}

\newcommand{\sgn}{\operatorname{sgn}}
\newcommand{\pfaff}{\mathop{\mathrm{pfaff}}}
\newcommand{\ep}{\varepsilon}

\numberwithin{equation}{section} 
\theoremstyle{plain}
\newtheorem{theo+}           {Theorem}      [section]
\newtheorem{prop+}  [theo+]  {Proposition}
\newtheorem{coro+}  [theo+]  {Corollary}
\newtheorem{lemm+}  [theo+]  {Lemma}
\newtheorem{defi+}  [theo+]  {Definition}

\theoremstyle{definition}
\newtheorem{exam+}  [theo+]  {Example}
\newtheorem{rema+}  [theo+]  {Remark}

\newenvironment{theorem}{\begin{theo+}}{\end{theo+}}
\newenvironment{proposition}{\begin{prop+}}{\end{prop+}}
\newenvironment{corollary}{\begin{coro+}}{\end{coro+}}
\newenvironment{lemma}{\begin{lemm+}}{\end{lemm+}}

\newenvironment{remark}{\begin{rema+}}{\end{rema+}}

\begin{document}

\baselineskip 18pt
\larger[2]

\title
[Sums of squares from elliptic pfaffians]
{Sums of squares
 from elliptic pfaffians} 
\author{Hjalmar Rosengren}
\address
{Department of Mathematics\\ Chalmers University of Technology and G\"oteborg
 University\\SE-412~96 G\"oteborg, Sweden}
\email{hjalmar@math.chalmers.se}
\urladdr{http://www.math.chalmers.se/{\textasciitilde}hjalmar}

\keywords{Representation of integers as sums of squares, pfaffian, Hankel determinant, theta function, elliptic function, Schur $Q$-polynomial,  Meixner--Pollaczek polynomial, continuous Hahn polynomial, correlation function, orthogonal polynomial ensemble}
\subjclass{Primary: 11E25, Secondary: 15A15, 33C45, 33E05}
\thanks{Research  supported by the Swedish Science Research
Council (Vetenskapsr\aa det)}

\begin{abstract}
We give a new proof of Milne's formulas for the number of representations of an integer as a sum of $4m^2$ and $4m(m+1)$ squares. The proof is based on explicit evaluation of pfaffians with elliptic function entries, and relates Milne's formulas to Schur $Q$-polynomials and to correlation functions for continuous dual Hahn polynomials. We  also state a new formula for  $2m^2$ squares.
\end{abstract}

\maketitle

 \section{Introduction}
One of the  classical problems of number theory is to count the number of representations of a  positive integer $n$ as a sum of $k$ squares. We will denote this number by $\square_k(n)$, where, as is customary,  representations
$$n=x_1^2+\dots+x_k^2 $$
 that may be obtained from each other by permuting the $x_i$, or replacing some $x_i$ by $-x_i$, are counted as different. 

The most fundamental results are Gauss' two squares and Jacobi's four and eight squares  formulas:
\begin{subequations}\label{s}
\begin{align}\label{s2}\square_2(n)&=4\sum_{d\mid n,\, d \text{ odd}}(-1)^{\frac12(d-1)}
,\\
\label{s4} \square_4(n)&=8\sum_{d\mid n,\,4\nmid d}d,\\
\label{s8} \square_8(n)&=
16\sum_{d\mid n}(-1)^{n+d} d^3,\end{align}
\end{subequations}
where the sums run over positive divisors. These may be compared with Legendre's formulas for sums of triangles:
\begin{subequations}\label{t}
\begin{align}\label{t2}\triangle_2(n)&= \sum_{d\mid
    4n+1}(-1)^{\frac 12(d-1)},\\
\label{t4}  \triangle_4(n)&=\sum_{d\mid 2n+1}d,\\
\label{t8}  \triangle_8(n)&=\sum_{ d\mid n+1,\,
  (n+1)/d \text{ odd}}d^3.
\end{align}
\end{subequations}

It is known that, in general,
$\square_{2k}(n)$ can be written as the sum of two terms, the first being an elementary divisor sum and the second the $n$-th Fourier coefficient of a  cusp form. Moreover, the second term vanishes only for $k\leq 4$, which indicates that there is no very simple extension of \eqref{s} to more than $8$ squares.

A novel approach, 
  motivated by  \emph{affine superalgebras}, was initiated by Kac and Wakimoto \cite{kw},
who used  \emph{denominator formulas} for such algebras to derive several new infinite families of  identities, both for squares and triangles. A particularly
 interesting case is  the ``queer'' series of affine superalgebras  $Q(m)$, for which the denominator formula was merely conjectured. For $Q(2m-1)$ and $Q(2m)$, respectively, it implies the triangular number identities 
\begin{subequations}\label{kmt}
\begin{equation}\label{kmt1}\triangle_{4m^2}(n)=\frac{1}{4^{m(m-1)}\prod_{j=1}^{2m-1}j!}
\sum_{\substack{k_1l_1+\dotsm+k_ml_m=2n+m^2\\k_1>k_2>\dots>k_m \\k_i \text{ and } l_i \text{ odd positive}}}
\,\prod_{i=1}^m k_i\prod_{1\leq
  i<j\leq m} (k_i^2-k_j^2)^2,\end{equation}
\begin{equation}\label{kmt2}\triangle_{4m(m+1)}(n)=
\frac{ 2^m}{\prod_{j=1}^{2m}j!}\sum_{\substack{k_1l_1+\dotsm
  +k_ml_m=n+\frac 12m(m+1)\\k_1>k_2>\dots>k_m \\ k_i\text{ positive},\ l_i \text{
  odd  positive}}}\,\prod_{i=1}^m k_i^3\prod_{1\leq
  i<j\leq m} (k_i^2-k_j^2)^2.\end{equation}
\end{subequations}
The case $m=1$ gives Legendre's formulas for four and eight 
triangles. The sums of squares formulas contained in \cite{kw} are different from those discussed in the present paper.

The identities \eqref{kmt} were first proved by Milne \cite{mi}, see also the research announcements \cite{mi1,mi2}, using an approach different from Kac and Wakimoto. Independently, the more general denominator formula was proved by Zagier \cite{z}.
  Getz and Mahlburg \cite{gm} showed that it contains 
 further triangular number identities, such as
 \begin{equation}\label{gm}\triangle_{2m}(n)=\sum_{\substack{k_1l_1+\dots+k_ml_m=4mn+m^2\\
     k_i \text{ {and} } l_i \text { {odd positive}}\\
     k_i\equiv\pm(2i-1)\!\mod 4m}}
 (-1)^{|\{i;\,k_i\equiv
   1-2i\ (4m)\}|},
  \end{equation}
 which reduces to \eqref{t2} for $m=1$.

In  \cite{r1}, we made the observation that
the denominator formulas for queer affine superalgebras can be written as \emph{pfaffian evaluations}, for matrices with elliptic function entries. 
Moreover, we showed that they follow  from a classical  determinant 
evaluation due to Frobenius \cite{f}. 
We also attempted  a  complete study of the implied triangular number identities,  finding formulas for $4m^2/d$ triangles, when $d\mid 2m$, and $4m(m+1)/d$ triangles, when $d\mid 2m$ or $d\mid 2m+2$. As an example, letting  $d=2$ in \cite[Eq.\ (4.4b)]{r1} gives a $2m^2$ triangles identity, which can be written as
\begin{multline}\label{hti}\triangle_{2m^2}(n)=\frac{(-1)^{\frac 12m(m-1)}}{4^{m(m-1)}\prod_{j=1}^{m-1}(j!)^2}\sum_{\substack{k_1l_1+\dotsm+k_ml_m=4n+m^2\\k_1>k_2>\dots>k_m \\k_i \text{ and } l_i \text{ odd positive}}}\prod_{i=1}^m(-1)^{\frac 12(k_i-1)}\\
\times\prod_{1\leq i<j\leq m}\left((-1)^{\frac12(k_j-1)}k_j-(-1)^{\frac 12(k_i-1)}k_i\right)^2.  \end{multline}


Milne's proof of \eqref{kmt} also uses elliptic functions, but is based on continued fractions and Hankel determinants, rather than on pfaffians. With similar methods he  obtained new formulas  for sums of $4m^2$ and $4m(m+1)$ squares, see Corollary \ref{mt} below. As an example, Milne's $16$ squares formula can be written \cite[Corollary~8.1]{mi}
\begin{multline*}
\square_{16}(n)=\frac{2^5}{3}\sum_{kl=n,\, k,l\geq 1}(-1)^{(k-1)(l-1)}k(1+k^2+k^4)\\
+\frac{2^8}{3}\sum_{\substack{k_1l_1+k_2l_2=n\\ k_1>k_2\geq 1,\,l_1,l_2\geq 1}}
(-1)^{(k_1-1)(l_1-1)+(k_2-1)(l_2-1)}k_1k_2(k_1^2-k_2^2)^2.
\end{multline*}
For a readable introduction to Milne's work, we refer to the  survey
 \cite{ck}. 

Ono \cite{o} derived seemingly different formulas for $4m^2$ and $4m(m+1)$ squares  from \eqref{kmt}, 
using the fact that the generating functions for squares and triangles are related by a modular transformation, cf.\  \eqref{mts} below. Moreover, in the introduction he indicates an alternative proof, also  using modular forms but without relying on \eqref{kmt}.
In the Appendix, we will show that Ono's formulas are in fact equivalent to Milne's, although his proof is  different. Another modular forms proof of Milne's formulas was recently given by Long and Yang \cite{ly}. This should be
 related to the alternative proof indicated by Ono.

The purpose of the present work is to extend the analysis of \cite{r1} from  triangles to squares, deriving Milne's sums of squares formulas from elliptic pfaffian evaluations. These evaluations follow from Frobenius' classical determinant, and may be viewed as modular duals of those used in \cite{r1}. Thus, just as in Ono's proof square number identities arise as modular duals of triangular number identities, although the details are  different. 

To be precise, we  give two, closely related, derivations of Milne's formulas. Applying two different modified Laurent series expansions to the matrix elements of our pfaffians, we obtain two different analogues of the Kac--Wakimoto denominator formula, Theorem \ref{adf} and Theorem \ref{mdt}. Theorem \ref{adf} looks very similar to the Kac--Wakimoto formula, though  convergent Lambert series are replaced by Abel means of divergent series. Theorem \ref{mdt} involves convergent series, but Schur polynomials are replaced by  \emph{Schur $Q$-polynomials}. To be precise, classical  Schur $Q$-polynomials are labelled by positive integer partitions; here, we need an extension to the case when some labels are negative. 

Theorems \ref{adf} and \ref{mdt} involve a number of free variables.  Specializing all  variables to~$1$ gives formulas for $4m^2$ and $4m(m+1)$ squares. The resulting formulas are equivalent, though to see that is far from obvious. Starting with Theorem \ref{adf}, one readily obtains the Hankel determinant form of Milne's identities, Corollary~\ref{mhd}. Applying an identity relating Hankel determinants to \emph{correlation functions} of orthogonal polynomial ensembles, we deduce explicit sums of squares formulas involving correlation functions for \emph{continuous dual Hahn polynomials}, see Corollary \ref{sst}. Expanding the correlation functions into Schur polynomials yields the  Schur function form of Milne's identities, Corollary \ref{mt}.

On the other hand, specializing all variables in Theorem \ref{mdt} to $1$ yields sums of squares formulas involving similarly specialized  Schur $Q$-polynomials, Corollary~\ref{qss}. 
Such quantities, and more generally Schur $Q$-polynomials specialized to a geometric progression, are studied in  \cite{r3}. It follows from that work that the Schur $Q$-polynomials appearing in Corollary \ref{qss} agree with the correlation functions of Corollary \ref{sst}, and consequently that Corollary \ref{qss} is equivalent to Milne's identities.

It is natural to try to generalize Milne's formulas to other numbers of squares, similarly as was done for triangular numbers in \cite{gm,r1}. 
In principle, it should be possible to obtain formulas for $4m^2/d$ squares, when $d\mid 2m$ and $4m(m+1)/d$ squares, when $d\mid 2m$ or $d\mid 2(m+1)$, from our elliptic pfaffian evaluations. However, in practice the computations are difficult to handle, and the resulting identities seem very complicated to state. 
As an indication of what the results may look like, we  state without proof a new formula for    $2m^2$ squares  in  Theorem \ref{hsf}. In particular, this embeds  the two squares formula \eqref{s2} in an infinite family, similarly as Milne's identities do for four and eight squares. 


The reader may question the merits of our new proof of Milne's formulas, compared to the  simple modular forms proofs of Ono and of  Long and Yang. First of all,  in our approach (similarly as in Milne's original proof) the identities are \emph{derived} from classical results. Although modular forms are a  powerful tool for verifying  Milne's formulas, they seem less useful for deriving the results from scratch.  
Second, our proof explains the  coefficients in Milne's formulas (i.e.\ the double sums over $\lambda$ and $\mu$ in Corollary \ref{mt}), by relating them to Schur $Q$-polynomials and to correlation functions for 
continuous dual Hahn polynomials.
Third, we obtain Milne's formulas as a special case of the  more general  Theorems \ref{adf} and \ref{mdt}, which might be of independent interest, plausibly in the theory of superalgebras. 
In conclusion, we feel that the three known approaches to Milne's identities (Milne's original proof using continued fractions and Hankel determinants, the modular forms proofs of   Ono and of  Long and Yang and the present work using pfaffian evaluations) each add something to the understanding of these deep and beautiful results.

{\bf Acknowledgements:} I thank Stephen Milne for several useful comments. 

\section{Preliminaries}

\subsection{Pfaffians}

As in \cite{r1}, we define the pfaffian of a  skew-symmetric matrix as
$$\pfaff_{1\leq i,j\leq m}(a_{ij})
=\frac{1}{2^{M} M!}\sum_{\sigma\in S_{m}}\sgn(\sigma)\prod_{i=1}^{M}
a_{\sigma(2i-1),\sigma(2i)},$$
where $M$ is the integral part of $m/2$. 
This definition is standard when $m$ is even, but not when $m$ is odd.
In even dimension,
\begin{equation}\label{pdi}\det_{1\leq i,j\leq 2m}(a_{ij})=\left(\pfaff_{1\leq i,j\leq 2m}(a_{ij})\right)^2. \end{equation}
The  cases of odd and even dimension are  related by
$$\pfaff_{1\leq i,j\leq 2m+1}(a_{ij})=\pfaff_{1\leq i,j\leq 2m+2}\left(\begin{matrix}a_{ij}&\begin{matrix}1\\\vdots \\ 1\end{matrix}
    \\\begin{matrix} -1 & \cdots & -1\end{matrix} & 0\end{matrix}\right).$$
We note in passing the pfaffian evaluation
\begin{equation}\label{spe}\pfaff_{1\leq i,j\leq m}\left(\frac{x_i-x_j}{x_i+x_j}\right)=\prod_{1\leq i<j\leq m}\frac{x_i-x_j}{x_i+x_j}. \end{equation}
When $m$ is even, this is a classical identity of Schur \cite{s}. The case of odd $m$ then follows by letting $x_m=0$.

\subsection{Theta functions}

Throughout, $q$ will be a fixed number with $0<q<1$. 
We will use the notation
$$(a)_\infty=(a;q)_\infty=\prod_{k=0}^\infty(1-aq^{k}),$$
$$(a_1,\dots,a_m)_\infty=(a_1,\dots,a_m;q)_\infty=(a_1;q)_\infty\dotsm(a_m;q)_\infty. $$
When the base  is suppressed from the notation, it is always taken to equal our fixed number $q$. 

We introduce the theta function
$$\theta(x)=\theta(x;q)=(x,q/x;q)_\infty.$$
It satisfies
\begin{equation}\label{qp}\theta(x^{-1})=\theta(qx)=-x^{-1}\theta(x) \end{equation}
and the modular transformation
\begin{multline}\label{mtt}\theta(e^{2\pi ix};e^{-2\pi/h})\\
=-i\sqrt
h\,e^{-\frac\pi{4}(h-h^{-1})}\frac{(e^{-2\pi h};e^{-2\pi h})_\infty}{(e^{-2\pi/h};e^{-2\pi/h})_\infty}\,e^{\pi x(i+h(1-x))}\,\theta(e^{-2\pi hx};e^{-2\pi h}).\end{multline}

The Laurent expansion of $\theta$ is given by Jacobi's triple product identity
\begin{equation}\label{jtp}(q)_\infty\,\theta(x)=\sum_{k=-\infty}^\infty(-1)^kq^{\binom k2}x^k.\end{equation}
We also mention the Laurent expansion 
\begin{equation}\label{rp}\frac{(q)_\infty^2\theta(ax)}{\theta(a)\theta(x)}=\sum_{k=-\infty}^\infty
\frac{x^k}{1-aq^k},\qquad q<|x|<1,
 \end{equation}
which is a special case of
 Ramanujan's ${}_1\psi_1$ sum \cite[Eq.\ (5.2.1)]{gr},
together with its limit case
\begin{equation}\label{ol}x\,\frac{\theta'(x)}{\theta(x)}=-\sum_{k\neq
  0}\frac{x^{k}}{1-q^k},\qquad q<|x|<1.\end{equation}

\subsection{Generating functions}

The triple product identity \eqref{jtp}
 implies explicit formulas for the generating functions for  squares and triangles. Indeed, let
$$\square(q)=\sum_{n=-\infty}^\infty q^{n^2}=1+2\sum_{n=1}^\infty q^{n^2}, $$
$$\triangle(q)=\frac 12\sum_{n=-\infty}^\infty q^{\frac12n(n+1)}=\sum_{n=0}^\infty q^{\frac12n(n+1)},$$
 so that, according to the standard conventions that we use,
$$
 \square(q)^k=\sum_{n=0}^\infty \square_k(n)q^n,$$
$$\triangle(q)^k=\sum_{n=0}^\infty \triangle_k(n)q^n. $$
Then, by \eqref{jtp},
\begin{equation}\label{sp}\square(q)=(q^2;q^2)_\infty\,\theta(-q;q^2)=(q^2,-q,-q;q^2)_\infty=\frac{(-q;-q)_\infty}{(q;-q)_\infty},\end{equation}
$$\triangle(q)=\frac 12\,(q;q)_\infty\,\theta(-q;q)=(q,-q,-q;q)_\infty=\frac{(q^2;q^2)_\infty}{(q;q^2)_\infty}.$$
Note that the special case $x=1/2+i/h$ of \eqref{mtt} gives
 \begin{equation}\label{mts}\triangle(e^{-2\pi/h})=\frac {\sqrt{h}}{2}\,e^{\pi/4h}\,\square(-e^{-\pi h}). \end{equation}

 We recall the Lambert series versions of \eqref{s}, that is,
 \begin{subequations}\label{l}
 \begin{equation}\label{l2}\square(q)^2=1+4\sum_{k=1}^\infty\frac{q^k}{1+q^{2k}},
 \end{equation}
 \begin{equation}\label{l4}\square(q)^4=1+8\sum_{k=1}^\infty\frac{kq^k}{1+(-q)^{k}}, \end{equation}
 \begin{equation}\label{l8}\square(q)^8=1+16\sum_{k=1}^\infty\frac{k^3q^k}{1-(-q)^{k}}. \end{equation}
 \end{subequations}
 Expanding the denominators as geometric series leads immediately to \eqref{s2} and \eqref{s8}. In the case of four squares, one obtains
  $$8\sum_{d\mid n} (-1)^{(d-1)(\frac nd-1)}d, $$
 which is  equivalent to \eqref{s4}; cf.\ the case   $k=1$ of \eqref{oe}  below.

In \S \ref{hss}, we  need  two further identities of Jacobi. 
First, we have the Lambert series expansion
 \begin{equation}\label{jl}q\triangle(q^2)^4\,\square(q)^2=\sum_{k=1}^\infty\frac{k^2q^k}{1+q^{2k}}. \end{equation}
 Second,  comparing \eqref{s4} and \eqref{t4} gives
 $$8\triangle_4(n)=\square_4(2n+1), $$
 which can equivalently be written
 \begin{equation}\label{jq}16q\triangle(q^2)^4=\square(q)^4-\square(-q)^4. \end{equation}

\subsection{Sums of triangular numbers}
\label{tss}

For  comparison, we briefly recall how the triangular number identities \eqref{kmt} follow from the pfaffian evaluations
\begin{subequations}\label{tp}
\begin{equation}\label{ep}
\pfaff_{1\leq i,j\leq 2m}\left(\frac{\theta(x_j/x_i)}{x_j\theta(\sqrt
    qx_j/x_i)}\right)= q^{\frac
    14m(m-1)}\prod_{i=1}^{2m}x_i^{m-i}\prod_{1\leq i<j\leq
    2m}\frac{\theta(x_j/x_i)}{\theta(\sqrt qx_j/x_i)},\end{equation}
\begin{multline}\label{op}
\pfaff_{1\leq i,j\leq 2m+1}\left(\frac{x_i\theta'(\sqrt
  qx_i/x_j)}{x_j\theta(\sqrt q x_i/x_j)}\right)\\
= q^{\frac
  14m(m-1)}\frac{(q)_\infty^{2m}}{(\sqrt
  q)_\infty^{2m}}\prod_{i=1}^{2m+1}x_i^{m+1-i}\prod_{1\leq i<j\leq
    2m+1}\frac{\theta(x_j/x_i)}{\theta(\sqrt qx_j/x_i)}.\end{multline}
\end{subequations}
For further details we refer to \cite{r1}, though the essential points are contained already in \cite{z}. 

It will be convenient to introduce the functions
\begin{equation}\label{sf}S_{\mu}(x_1,\dots,x_m)=\frac{\det_{1\leq i,j\leq m}(x_j^{\mu_i})}{\prod_{1\leq i<j\leq m}(x_i-x_j)}, \qquad \mu\in\mathbb Z^m, \end{equation}
which are essentially Schur polynomials. Indeed, since $S_\mu$ is anti-symmetric in the variables $\mu_i$ we may assume $\mu_1>\dots>\mu_m$. Moreover, since
\begin{equation}\label{st}S_{(\mu_1+a,\dots,\mu_m+a)}(x_1,\dots,x_m)=x_1^a\dotsm x_m^a\, S_{(\mu_1,\dots,\mu_m)}(x_1,\dots,x_m) \end{equation}
we may also assume $\mu_m\geq 0$. Then, writing $\mu_i=\lambda_i+m-i$, so that $\lambda_1\geq\dots\geq\lambda_m\geq 0$,  $S_\mu$ equals the Schur polynomial $s_\lambda$.

Returning to the task at hand, we
 expand each matrix element in \eqref{tp} as a Laurent series in the annulus $\sqrt q<|x_j/x_i|<1/\sqrt q$. Applying \eqref{rp} and \eqref{ol} one arrives 
after some elementary manipulation at the multivariable Lambert series
\begin{subequations}\label{dfe}
\begin{multline}
\frac{(q)_\infty^{2m}}{(\sqrt q)_\infty^{2m}}\prod_{i=1}^{2m}\frac 1{x_i^m}\prod_{1\leq i<j\leq 2m}\frac{(qx_j/x_i,qx_i/x_j)_\infty}{(\sqrt q x_j/x_i,\sqrt q x_i/x_j)_\infty}
=q^{-\frac 14m(m-1)}\\
\times\sum_{k_1>\dots>k_m\geq 0}\,\prod_{i=1}^m\frac{q^{\frac 12 k_i}}{1-q^{k_i+\frac 12}}\,S_{(k_1,\dots,k_m,-k_m-1,\dots,-k_1-1)}(x_1,\dots,x_{2m}),
\end{multline}
\begin{multline}
\frac{(q)_\infty^{2m}}{(\sqrt q)_\infty^{2m}}\prod_{i=1}^{2m+1}\frac 1{x_i^m}\prod_{1\leq i<j\leq 2m+1}\frac{(qx_j/x_i,qx_i/x_j)_\infty}{(\sqrt q x_j/x_i,\sqrt q x_i/x_j)_\infty}\\
=q^{-\frac 14m(m+1)}\sum_{k_1>\dots>k_m\geq 1}\,\prod_{i=1}^m\frac{q^{\frac 12 k_i}}{1-q^{k_i}}\,S_{(k_1,\dots,k_m,0,-k_m,\dots,-k_1)}(x_1,\dots,x_{2m+1}),
\end{multline}
\end{subequations}
which are  equivalent to the denominator formula for $Q(2m-1)$ and $Q(2m)$, respectively. 

Specializing $x_i\equiv 1$, the left-hand sides of \eqref{dfe}
reduce to $\triangle(\sqrt q)^{4m^2}$ and $\triangle(\sqrt q)^{4m(m+1)}$, respectively. On the right, applying the classical formula 
\begin{equation}\label{hlf}S_\mu(1^m)=\prod_{1\leq i<j\leq m}\frac{\mu_i-\mu_j}{j-i}\end{equation}
  and expanding the denominators as geometric series leads after simplifications to \eqref{kmt}.

More generally, one may let $x_i=\omega^{i-1}$, with  $\omega$ a suitable root of unity. This leads to  more general triangular number identities  such as \eqref{gm} and \eqref{hti}, see \cite{r1}. 

\subsection{Hankel determinants}

We recall the following classical result  \cite[Corollary 2.1.3]{i}.

\begin{lemma}\label{hdl}
Let
$$\mu(f)=\int f(x)\,d\mu(x)$$
be a linear functional defined on polynomials, let $c_k=\mu(x^k)$ be its moments, and let
$$\Delta(x)=\prod_{1\leq i<j\leq m}(x_j-x_i)=\det_{1\leq i,j\leq m}\left(x_i^{j-1}\right)
$$
denote the Vandermonde determinant. Then,
$$\det_{1\leq i,j\leq m}\left(c_{i+j-2}\right)
=
\frac 1{m!}\int \Delta(x_1,\dots,x_m)^2\,d\mu(x_1)\dotsm d\mu(x_m).$$
\end{lemma}

\begin{proof}
Clearly,
$$\int \Delta(x_1,\dots,x_m)^2\,d\mu(x_1)\dotsm d\mu(x_m)=\sum_{\sigma,\tau\in S_m}\sgn(\sigma)\sgn(\tau)\prod_{i=1}^mc_{\sigma(i)+\tau(i)-2}.$$
Replacing $\sigma$ by $\sigma\tau$, this may indeed be written
$$\sum_{\sigma,\tau\in S_m}\sgn(\sigma)\prod_{i=1}^m c_{i+\sigma(i)-2}=m!\det_{1\leq i,j\leq m}(c_{i+j-2}). $$
\end{proof}

Equivalently, defining the  coefficients $C(k_1,\dots,k_m)$ by
$$\prod_{1\leq i<j\leq m}(x_j-x_i)^2=\sum_{k_1,\dots,k_m}C(k_1,\dots,k_m)\prod_{i=1}^m x_i^{k_i}, $$
we have, for arbitrary scalars $c_k$,
$$\sum_{k_1,\dots,k_m}C(k_1,\dots,k_m)\prod_{i=1}^m c_{k_i}=m!\det_{1\leq i,j\leq m}(c_{i+j-2}). $$
This version of Lemma \ref{hdl} will be useful for our discussion of Ono's identities in the Appendix.

If $p_k$ and $q_k$ are arbitrary monic polynomials of degree $k$ then, by linearity,
$$\det_{1\leq i,j\leq m}\left(\mu(p_{i-1}q_{j-1})\right)=\det_{1\leq i,j\leq m}\left(\mu(x^{i+j-2})\right). $$
In particular, if $\mu$ is a positive functional we may choose $p_k=q_k$ as the corresponding monic orthogonal polynomials. Then, Lemma \ref{hdl} gives
\begin{equation}\label{hop}
\frac 1{m!}\int \Delta(x_1,\dots,x_m)^2\,d\mu(x_1)\dotsm d\mu(x_m)
=\prod_{i=1}^m \|p_{i-1}\|^2.\end{equation}

\subsection{Tangent numbers}
We will need the following evaluation of  Abel means of alternating power sums in terms of  tangent numbers, see \cite[Theorem 2.5]{t}. It is  equivalent to the classical evaluation of Riemann's zeta function at the negative integers, though for completeness we provide a self-contained proof.
We find it convenient to define  $t_k$ by
\begin{equation}\label{cb}\tan\frac x2=\sum_{k=1}^\infty t_k\frac{ x^{2k-1}}{(2k-1)!}. \end{equation}
Equivalently, in standard notation for tangent and Bernoulli numbers,
$$t_k=2^{1-2k}T_k=\frac{(4^{k}-1)|B_{2k}|}{k}.$$
We will denote Abel means of possibly divergent series by
\begin{equation}\label{am}\sideset{}{'}\sum_{k\in\Lambda}c_k=\lim_{t\rightarrow 1^-}\sum_{k\in\Lambda}t^{|k|}c_k,\qquad \Lambda\subseteq \mathbb Z. \end{equation}

\begin{lemma}
One has
$$\sideset{}{'}\sum_{k=1}^\infty(-1)^{k+1} k^m=\begin{cases}
1/4, & m=0,\\
0, & m=2,4,6,\dots,\\
(-1)^{n+1}t_n/2, & m=2n-1.
\end{cases} $$
\end{lemma}

\begin{proof}
We assume that $m>0$. Since
 \begin{equation}\label{rfe}\frac{1-x}{1+x}=1+2\sum_{k=1}^\infty(-1)^kx^k,
 \qquad |x|<1,\end{equation}
we may write
\begin{equation*}\begin{split}2\sideset{}{'}\sum_{k=1}^\infty(-1)^{k+1} k^m&=\left(x\frac{d}{dx}\right)^m
\Bigg|_{x=1}\frac{1-x}{1+x}=\left(\frac{d}{dt}\right)^m\Bigg|_{t=0}
\frac{1-e^t}{1+e^t}\\
&=\left(\frac{d}{dt}\right)^m\Bigg|_{t=0}\left(
i\tan\frac t{2i}
\right),\end{split}\end{equation*}
which implies the desired result.
\end{proof}

We introduce the moment functionals
\begin{equation}\label{mue}\mu_\ep(f)= 2\sideset{}{'}\sum_{k=1}^\infty(-1)^{k+1+\ep}  k^{1+2\ep}\,f(-k^2),\qquad \ep=0,\,1.  \end{equation}

\begin{corollary}\label{muc}
The moments of $\mu_0$ and $\mu_1$ are given by
\begin{equation}\label{muce}\mu_0(x^k)=\mu_1(x^{k-1})=t_{k+1}.\end{equation}
\end{corollary}

 In Lemma \ref{tml},
we  give alternative expressions for the functionals $\mu_\ep$. As a consequence, we shall see that  the corresponding orthogonal polynomials are continuous dual Hahn polynomials.  

\section{From pfaffians to Hankel determinants}
\label{fps}

We will obtain sums of squares formulas from the following pfaffian evaluations, the first of which was already given in \cite{r1}.

\begin{lemma}\label{epe}
One has
\begin{subequations}\label{epi}
\begin{equation}\label{eep}\pfaff_{1\leq i,j\leq 2m}\left(\frac{\theta(x_j/x_i)}{\theta(-x_j/x_i)}\right)=\prod_{1\leq i<j\leq
    2m}\frac{\theta(x_j/x_i)}{\theta(-x_j/x_i)},\end{equation}
\begin{equation}\label{oep}\pfaff_{1\leq i,j\leq 2m+1}\left(1+2\frac{x_j\,\theta'(-x_j/x_i)}{x_i\,\theta(-x_j/x_i)}\right)=\frac{(q)_\infty^{2m}}{(-q)_\infty^{2m}}\prod_{1\leq i<j\leq
    2m+1}\frac{\theta(x_j/x_i)}{\theta(-x_j/x_i)}.\end{equation}
\end{subequations}
\end{lemma}

\begin{remark}\label{r2}
Since $\theta(x;0)=1-x$,  Lemma \ref{epe} reduces to  \eqref{spe} when $q=0$. A different elliptic extension of Schur's pfaffian evaluation was recently obtained by Okada \cite{ok}, see also 
\cite[Remark 2.1]{r1}. 
\end{remark}

We  give two proofs of \eqref{oep}. The identity \eqref{eep} can be proved similarly \cite{r1}.

\begin{proof}[First proof of \eqref{oep}]
 We start from the Frobenius--Stickelberger determinant \cite{fs}
\begin{multline*} \det_{1\leq i,j\leq
  n+1}\left(\begin{matrix}\displaystyle-\frac{x_jy_i\,\theta'(x_jy_i)}{(q)_\infty^2\theta(x_jy_i)}&\begin{matrix}1\\
  \vdots \\ 1\end{matrix}\\ \begin{matrix}
-1 & \cdots & -1 \end{matrix} & 0\end{matrix}\right)\\
  =\frac{\theta(x_1\dotsm x_ny_1\dotsm y_n)\prod_{1\leq i<j\leq
  n}x_jy_j\theta(x_i/x_j)\theta(y_i/y_j)}{\prod_{i,j=1}^n\theta(x_iy_j)}.
\end{multline*}
Choosing $n=2m+1$ and $y_j=-1/x_j$, and moreover adding $1/2(q)_\infty^2$ times the last column to all other columns gives
\begin{multline*}\det_{1\leq i,j\leq
  2m+2}\left(\begin{matrix}\displaystyle \frac{1}{2(q)_\infty^2}
\left(1+2\frac{x_j\,\theta'(-x_j/x_i)}{x_i\,\theta(-x_j/x_i)}\right)
&\begin{matrix}1\\
  \vdots \\ 1\end{matrix}\\ \begin{matrix}
-1 & \quad\qquad\cdots\quad\qquad & -1 \end{matrix} & 0\end{matrix}\right)\\
  =\frac{1}{4^{m}(-q)_\infty^{4m}}\prod_{1\leq i<j\leq
  n}\frac{\theta(x_i/x_i)^2}{\theta(-x_j/x_i)^2}.
\end{multline*}
We now observe that the matrix on the left is skew-symmetric. 
This is easily proved by differentiating \eqref{qp}. 
  Thus, in view of
\eqref{pdi}, we may conclude that \eqref{oep} holds up to a factor $\pm 1$, which is independent of $q$ by continuity. By Remark \ref{r2},  that factor has to be $+1$.
\end{proof}

\begin{proof}[Second proof of \eqref{oep}]
In \eqref{op}, let $q=e^{-2\pi/h}$ and $x_j=e^{2\pi iz_j}$. 
Applying the logarithmic derivative of \eqref{mtt}, that is,
$$2ie^{2\pi ix}\frac{\theta'(e^{2\pi ix};e^{-2\pi/h})}{\theta(e^{2\pi ix};e^{-2\pi/h})}=i+h-2hx-2he^{-2\pi hx}\frac{\theta'(e^{-2\pi hx};e^{-2\pi h})}{\theta(e^{-2\pi hx};e^{-2\pi h})}, $$
the left-hand side of \eqref{op} takes the form
$$\left(\frac{h e^{\pi/h}}{2i}\right)^m\pfaff_{1\leq i,j\leq 2m+1}\left(
1+2\frac{e^{2\pi hz_j}\,\theta'(-e^{2\pi h(z_j-z_i)};e^{-2\pi h})}{e^{2\pi hz_i}\,\theta(-e^{2\pi h(z_j-z_i)};e^{-2\pi h})}+2(z_j-z_i)\right). $$
By Lemma \ref{rol} below, we may subtract $2(z_j-z_i)$ from each matrix element, which gives a pfaffian of the desired type. Similarly applying \eqref{mtt} and \eqref{mts} to the right-hand side of \eqref{op}  completes the proof.
\end{proof}

The  following lemma was used above.

\begin{lemma}\label{rol}
For any odd-dimensional skew-symmetric matrix $(a_{ij})$,
$$\pfaff_{1\leq i,j\leq 2m+1}(a_{ij}+b_i-b_j)=\pfaff_{1\leq i,j\leq 2m+1}(a_{ij}). $$
\end{lemma}

\begin{proof}
Adding  multiples of the last row and column to the previous ones gives
$$
\det_{1\leq i,j\leq 2m+2}\left(\begin{matrix}a_{ij}+b_i-b_j&\begin{matrix}1\\\vdots \\ 1\end{matrix}
    \\\begin{matrix} -1 & \cdots & -1\end{matrix} & 0\end{matrix}\right)=
\det_{1\leq i,j\leq 2m+2}\left(\begin{matrix}a_{ij}&\begin{matrix}1\\\vdots \\ 1\end{matrix}
    \\\begin{matrix} -1 & \cdots & -1\end{matrix} & 0\end{matrix}\right).$$
By \eqref{pdi}, we may extract square roots to conclude that
$$\pfaff_{1\leq i,j\leq 2m+1}(a_{ij}+b_i-b_j)=\pm\pfaff_{1\leq i,j\leq 2m+1}(a_{ij}).$$
Since both sides are polynomial expressions, the sign may be determined by letting $b_i\equiv 0$.
\end{proof}

Proceeding in analogy with \S \ref{tss}, the next step would be to expand the matrix elements in \eqref{epi} as Laurent series in the variables $x_i$. 
An important difference from \eqref{tp} is the presence of singularities at $x_i=-x_j$, which precludes Laurent expansion near $x_i=x_j=1$. To circumvent this difficulty we will use Abel means, for which we recall the notation \eqref{am}.

\begin{lemma}\label{asl} 
If $|x|=1$ and $x\neq 1$, then
\begin{subequations}\label{ata}
\begin{equation}\label{atea}\frac{(q)_\infty^2\theta(x)}{(-q)_\infty^2\theta(-x)}=\sideset{}{'}\sum_{k\neq 0}
 \frac{1-q^k}{1+q^k}\,(-x)^k,\end{equation}
\begin{equation}\label{atoa}1+2x\frac{\theta'(-x)}{\theta(-x)}=
\sideset{}{'}\sum_{k\neq 0}
 \frac{1+q^k}{1-q^k}\,(-x)^k. \end{equation}
\end{subequations}
\end{lemma}

\begin{proof}
By \eqref{rp},
\begin{equation}\label{rpa}\frac{(q)_\infty^2\theta(x)}{(-q)_\infty^2\theta(-x)}=2\sum_{k=-\infty}^\infty
\frac{(-x)^k}{1+q^k},\qquad q<|x|<1, \end{equation}
Replacing $x$ by $qx$  and using \eqref{qp} gives
\begin{equation}\label{rpb}\frac{(q)_\infty^2\theta(x)}{(-q)_\infty^2\theta(-x)}=-2\sum_{k=-\infty}^\infty
\frac{(-qx)^k}{1+q^k},\qquad 1<|x|<q^{-1}.\end{equation}
Letting $f(x)$ denote the left-hand side of \eqref{atea}, we write
$$f(x)=\lim_{t\rightarrow 1^-}\frac{f(xt)+f(x/t)}{2}. $$
Expanding $f(xt)$ using \eqref{rpa} and $f(x/t)$ using \eqref{rpb} yields \eqref{atea}. 

As for \eqref{atoa}, if
 $f(x)$ denote its left-hand side then 
\eqref{ol} gives
$$f(x)=1+2\sum_{k\neq 0}\frac{(-x)^k}{1-q^k},\qquad q<|x|<1. $$
It is easy to check that $f(1/x)=-f(x)$. Indeed, this is equivalent to the skew-symmetry of the matrix \eqref{oep}. Thus,
$$f(x)=-1+2\sum_{k\neq 0}\frac{q^k(-x)^k}{1-q^k},\qquad 1<|x|<q^{-1}. $$
Similarly as above, we can now deduce \eqref{atoa}.
\end{proof}

Applying \eqref{atea} to the matrix elements of
  \eqref{eep} gives
\begin{multline*}\frac{(q)_\infty^{2m}}{(-q)_\infty^{2m}}\prod_{1\leq i<j\leq
    2m}\frac{\theta(x_j/x_i)}{\theta(-x_j/x_i)}\\
\begin{split}&=\frac 1{2^mm!}\sum_{\sigma\in S_{2m}}\sgn(\sigma)\prod_{i=1}^m\left(\sideset{}{'}\sum_{k\neq 0}(-1)^k \frac{1-q^k}{1+q^k}\left(\frac{x_{\sigma(2i)}}{x_{\sigma(2i-1)}}\right)^k\right)\\
&=\frac 1{2^mm!}\sideset{}{'}\sum_{k_1,\dots,k_m\neq 0}\,\prod_{i=1}^m(-1)^{k_i}\frac{1-q^{k_i}}{1+q^{k_i}}\sum_{\sigma\in S_{2m}}\sgn(\sigma)\prod_{i=1}^m\left(\frac{x_{\sigma(2i)}}{x_{\sigma(2i-1)}}\right)^{k_i},\end{split}\end{multline*}
where we assume that  $|x_i|=|x_j|$ and $x_i\neq -x_j$ for all $i$ and $j$.
Multiplying both sides with $\prod_{i<j}(x_i+x_j)/(x_i-x_j)$ and adopting the notation \eqref{sf}, this may be written
\begin{multline*}\frac{(q)_\infty^{2m}}{(-q)_\infty^{2m}}\prod_{1\leq i<j\leq
    2m}\frac{(qx_j/x_i,qx_i/x_j)_\infty}{(-qx_j/x_i,-qx_i/x_j)_\infty}
=\frac 1{2^mm!}\prod_{1\leq i<j\leq 2m}(x_i+x_j)\\
\times\sideset{}{'}\sum_{k_1,\dots,k_m\neq 0}\prod_{i=1}^m(-1)^{k_i}\frac{1-q^{k_i}}{1+q^{k_i}}S_{(-k_1,k_1,\dots,-k_m,k_m)}(x_1,\dots,x_{2m}).
\end{multline*}
Since the summand is anti-symmetric and even as a function of $k_i$, 
and vanishes unless the $k_i$ are all distinct,
we may restrict the summation to $k_1>\dots>k_m\geq 1$ if we multiply by $2^m m!$. Moreover, 
$$
S_{(-k_1,k_1,\dots,-k_m,k_m)}=(-1)^mS_{(k_1,\dots,k_m,-k_m,\dots,-k_1)}.$$
This gives the first half of following result, which should be compared both with the Kac--Wakimoto identities \eqref{dfe} and with the closely related expansion given in
Theorem \ref{mdt} below. The second half follows similarly from 
\eqref{oep}.

\begin{theorem}\label{adf}
If $|x_i|=|x_j|$ and $x_i\neq -x_j$ for all $i$ and $j$, then
\begin{subequations}\label{ade}
\begin{multline}\label{adea}\frac{(q)_\infty^{2m}}{(-q)_\infty^{2m}}\prod_{1\leq i<j\leq
    2m}\frac{(qx_j/x_i,qx_i/x_j)_\infty}{(-qx_j/x_i,-qx_i/x_j)_\infty}
=\prod_{1\leq i<j\leq 2m}(x_i+x_j)\\
\times\lim_{t \rightarrow 1^-}\sum_{k_1>\dots>k_m\geq 1}\prod_{i=1}^mt^{k_i}(-1)^{k_i+1}\frac{1-q^{k_i}}{1+q^{k_i}}\,S_{(k_1,\dots,k_m,-k_m,\dots,-k_1)}(x_1,\dots,x_{2m}),
\end{multline}
\begin{multline}\label{adeb}\frac{(q)_\infty^{2m}}{(-q)_\infty^{2m}}\prod_{1\leq i<j\leq
    2m+1}\frac{(qx_j/x_i,qx_i/x_j)_\infty}{(-qx_j/x_i,-qx_i/x_j)_\infty}
=\prod_{1\leq i<j\leq 2m+1}(x_i+x_j)\\
\times\lim_{t \rightarrow 1^-}\sum_{k_1>\dots>k_m\geq 1}\prod_{i=1}^mt^{k_i}(-1)^{k_i}\frac{1+q^{k_i}}{1-q^{k_i}}\,S_{(k_1,\dots,k_m,0,-k_m,\dots,-k_1)}(x_1,\dots,x_{2m+1}).
\end{multline}
\end{subequations}
\end{theorem}

To obtain sums of squares formulas from  \eqref{ade}, we let $x_i\equiv 1$.  
 Then, the left-hand sides  reduce to $\square(-q)^{4m^2}$ and  $\square(-q)^{4m(m+1)}$. On the right, \eqref{hlf} gives
$$S_{(k_1,\dots,k_m,-k_m,\dots,-k_1)}(1^{2m})=\frac{2^m}{\prod_{i=1}^{2m-1}i!}
\prod_{i=1}^m k_i\prod_{1\leq i<j\leq m}(k_j^2-k_i^2)^2,
 $$
$$S_{(k_1,\dots,k_m,0,-k_m,\dots,-k_1)}(1^{2m+1})=\frac{2^m}{\prod_{i=1}^{2m}i!}
\prod_{i=1}^m k_i^3\prod_{1\leq i<j\leq m}(k_j^2-k_i^2)^2 .$$

Focusing on \eqref{adea},  we thus obtain
$$\square(q)^{4m^2}=\frac{4^{m^2}}{m!\prod_{i=1}^{2m-1}i!}
\sideset{}{'}\sum_{k_1,\dots,k_m= 1}^\infty\prod_{i=1}^m(-1)^{k_i+1}\frac{1-(-q)^{k_i}}{1+(-q)^{k_i}}\, k_i\prod_{1\leq i<j\leq m}(k_j^2-k_i^2)^2.$$
By Lemma \ref{hdl}, this may be written in Hankel determinant form as
$$\square(q)^{4m^2}=\frac{2^{m(2m-1)}}{\prod_{i=1}^{2m-1}i!}\det_{1\leq i,j\leq m}\left(\nu_0(x^{i+j-2})\right), $$
where
$$\nu_0(f)=2\sideset{}{'}\sum_{k=1}^\infty\frac{1-(-q)^k}{1+(-q)^k}\,(-1)^{k+1}kf(k^2). $$

Next, writing
$$\frac{1-(-q)^k}{1+(-q)^k}=1-2\frac{(-q)^k}{1+(-q)^k},
$$
 leads to the decomposition $\nu_0=\lambda_0+\mu_0$,  where $\mu_0$ is as in \eqref{mue}, and 
$$\lambda_0(f)=4\sum_{k=1}^\infty \frac{q^{k}k}{1+(-q)^{k}}\,f(-k^2);$$
since  the latter sum is convergent, there is no need to take the Abel mean. Using Corollary \ref{muc} to identify the moments of $\mu_0$, we arrive at the first half of   Corollary~\ref{mhd}.
 
Similarly, \eqref{adeb} leads to the identity
$$\square(q)^{4m(m+1)}=\frac{2^{m(2m+1)}}{\prod_{i=1}^{2m}i!}\det_{1\leq i,j\leq m}\left(\nu_1(x^{i+j-2})\right), $$
where
$$\nu_1(f)=2\sideset{}{'}\sum_{k=1}^\infty\frac{1+(-q)^k}{1-(-q)^k}\,(-1)^{k}k^3f(k^2). $$
Writing $\nu_1=\lambda_1+\mu_1$, with
 $\mu_1$  as in \eqref{mue} and
$$\lambda_1(f)=4\sum_{k=1}^\infty \frac{q^{k}k^3}{1-(-q)^{k}}\,f(-k^2).$$
 yields the second half of Corollary \ref{mhd}.  These identities are equivalent to  \cite[Theorems 5.3 and~5.5]{mi}.

\begin{corollary}[Milne]\label{mhd}
One has
$$\square(q)^{4m^2}=\frac{2^{m(2m-1)}}{\prod_{i=1}^{2m-1}i!}\,\det_{1\leq i,j\leq m}\left(t_{i+j-1}+4(-1)^{i+j}\sum_{k=1}^\infty \frac{q^{k}k^{2i+2j-3}}{1+(-q)^{k}}
\right),$$
$$\square(q)^{4m(m+1)}=\frac{2^{m(2m+1)}}{\prod_{i=1}^{2m}i!}\,\det_{1\leq i,j\leq m}\left(t_{i+j}+4(-1)^{i+j}\sum_{k=1}^\infty \frac{q^{k}k^{2i+2j-1}}{1-(-q)^{k}}
\right),$$
where the numbers $t_k$ are defined in \eqref{cb}.
\end{corollary}

Let $p_k^{(\ep)}(x;q)$ be the monic orthogonal polynomials corresponding to $\nu_\ep$. Then, \eqref{hop}
 gives the following reformulation of  Milne's  formulas. 

\begin{corollary}\label{opc}
In the notation above,
$$\square(q)^{4m^2}=\frac{2^{m(2m-1)}}{\prod_{i=1}^{2m-1}i!}\,\prod_{i=1}^m\|p_{i-1}^{(0)}(x;q)\|^2,
$$
$$\square(q)^{4m(m+1)}=\frac{2^{m(2m+1)}}{\prod_{i=1}^{2m}i!}\,\prod_{i=1}^m\|p_{i-1}^{(1) }(x;q)\|^2.$$
\end{corollary}

Equivalently,
\begin{subequations}\label{pqn}
\begin{equation}\|p_k^{(0)}(x;q)\|^2=\frac{(2k+1)!(2k)!}{2^{4k+1}}\,\square(q)^{8k+4},
 \end{equation}
\begin{equation}\|p_k^{(1)}(x;q)\|^2=\frac{(2k+2)!(2k+1)!}{2^{4k+3}}\,\square(q)^{8k+8}.\end{equation}
\end{subequations}

It is well-known that  any positive definite Hankel determinant can be expressed in terms of orthogonal polynomials. The point here is the explicit expressions for the moment functionals, which might suggest that
the  polynomials $p_k^{(\ep)}(x;q)$ are of  independent interest. In particular, an alternative proof of  \eqref{pqn} would lead to a new proof of Milne's sums of squares formulas. Note also that  $p_k^{(\ep)}(x;q)$ are $q$-analogues of the polynomials  $p_k^{(\ep)}(x;0)$, which, as we shall
 see in \S \ref{cdhs},  
are continuous dual Hahn polynomials. However, they are  of a different type from the polynomials in the $q$-Askey Scheme \cite{ks}.

\section{Correlation functions and sums of squares}

To obtain sums of squares formulas from  Corollary 
\ref{mhd}, we must expand the determinants as power series in $q$. In \cite{mi}, this is achieved using the Cauchy--Binet formula, which leads to  identities involving Schur functions, see Corollary~\ref{mt}. We will give a more conceptual derivation of these identities by relating them to   correlation functions for continuous dual Hahn polynomials.

\subsection{Continuous dual Hahn polynomials}
\label{cdhs}

We need the following alternative expressions for the functionals $\mu_\ep$ defined in \eqref{mue}.

\begin{lemma}\label{tml}
One has
 $$\mu_{\ep}(f)=\frac 12\int_{0}^\infty \frac{x^{2\ep+1} f(x^2)}{\sinh(\pi  x)}\,dx,\qquad  \ep=0,\,1.$$
\end{lemma}

\begin{proof}
Write
$$\frac 12\int_{0}^\infty \frac{x^{2\ep+1} f(x^2)}{\sinh(\pi  x)}\,dx
=\lim_{\lambda\rightarrow 0^+}\int_{-\infty}^\infty \frac{x^{2\ep+1}e^{i\lambda x} f(x^2)}{\sinh(\pi  x)}\,dx
 $$
and  expand  as the sum of  residues in the upper half-plane. 
\end{proof}

Alternatively, Lemma \ref{tml} can be deduced from Corollary \ref{muc}. 
One is then reduced  to the integral evaluation
$$\int_0^\infty\frac{x^{2k-1}}{\sinh(\pi x)}\,dx=\frac{t_k}{2},
$$
which can be found in standard tables, or  derived from 
$$\int_0^\infty\frac{\sinh(tx)}{\sinh(\pi x)}\,dx=\frac{\tan(t/2)}{2},\qquad -\pi<t<\pi $$
 by expanding both sides as power series in $t$.

Let  $(p_k^{(\ep)})_{k=0}^\infty$ be the monic orthogonal polynomials corresponding to $\mu_\ep$, $\ep=0,1$. We claim that they may be   identified  with  \emph{continuous dual Hahn polynomials} \cite{ks}. 
In general, when $a,b,c\geq 0$, one denotes by $(-1)^kS_k(x;a,b,c)$ the monic orthogonal polynomials with respect to 
$$f\mapsto\int_0^\infty \left|\frac{\Gamma(a+ix)\Gamma(b+ix)\Gamma(c+ix)}{\Gamma(2ix)} \right|^2\,f(x^2)\,dx. $$
Using elementary properties of the gamma function, it follows from Lemma \ref{tml} that
$$p_k^{(0)}(x)=(-1)^k S_k(x;0,1/2,1),\qquad 
p_k^{(1)}(x)=(-1)^k S_k(x;1/2,1,1).$$

Since $p_k^{(\ep)}(x)=p_k^{(\ep)}(x;0)$, 
 the case $q=0$ of Corollary \ref{opc}  reads 
\begin{equation}\label{en}
\prod_{i=1}^m\|p_{i-1}^{(0)}\|^2=\frac 1{2^{m(2m-1)}}\prod_{i=1}^{2m-1}i!,\qquad
\prod_{i=1}^m\|p_{i-1}^{(1)}\|^2=\frac 1{2^{m(2m+1)}}\prod_{i=1}^{2m}i!,
 \end{equation}
which agrees with the known expressions for the norms  \cite{ks}. 

We remark that
 $p_k^{(0)}$ and $p_k^{(1)}$    combine naturally to a single orthogonal system. Indeed, if 
$$p_{2k}(x)=p_k^{(0)}(x^2),\qquad p_{2k+1}(x)=xp_k^{(1)}(x^2),$$
then $(p_k(x))_{k=0}^\infty$ are the monic orthogonal polynomials corresponding to 
$$\int_{-\infty}^\infty \frac{xf(x)}{\sinh(\pi x)}\,dx. $$ 
These are Meixner--Pollaczek polynomials; in the notation of \cite{ks}
$$p_k(x)=\frac{k!}{2^k}\,P^{(1)}_k(x;\pi/2).$$

\subsection{Correlation functions}
\label{cfs}

We recall some general facts about correlation functions, see \cite{r2} for references and further details. 

Let 
$$\mu(f)=\int f(x)\,d\mu(x)$$
be a positive moment functional and $(p_k(x))_{k=0}^\infty$  the corresponding family of monic orthogonal polynomials. Normalizing 
$$\Delta(x_1,\dots,x_n)^2\,d\mu(x_1)\dotsm d\mu(x_m)$$
to a probability measure, it
defines an \emph{orthogonal polynomial ensemble}. 
Such ensembles arise in several contexts, including the theory of random hermitian matrices \cite{ko}. An important role is played by the correlation functions
$$\int\Delta(x_1,\dots,x_n)^2\,d\mu(x_{m+1})\dotsm d\mu(x_{n}),\qquad 0\leq m\leq n. $$
We will  normalize these functions as
\begin{multline}\label{kcf}C(x)= C_m^n(x_1,\dots,x_m)
=\frac{1}{(n-m)!\prod_{i=1}^n\| p_{i-1}\|^2\Delta(x_1,\dots,x_m)^2}\\ 
\times\int\Delta(x_1,\dots,x_n)^2\,d\mu(x_{m+1})\dotsm d\mu(x_{n}).
\end{multline}
Then, $C$ is a symmetric polynomial.

There are several expressions for correlation functions in terms of orthogonal polynomials, including
\begin{subequations}\label{cp}
\begin{gather}
\label{cd}C(x)=\frac{1}{\|p_{n-1}\|^{2m}\Delta(x)^2}\,\det_{1\leq i,j\leq m}\left(\frac{p_n(x_i)p_{n-1}(x_j)-p_{n-1}(x_i)p_n(x_j)}{x_i-x_j}\right)\\
 =\frac{(-1)^{\frac 12m(m-1)}}{\prod_{i=1}^m\|p_{n-i}\|^2\Delta(x)^4}\,\det_{1\leq i,j\leq 2m}\left(\begin{cases}p_{n-m+j-1}(x_i), & 1\leq i\leq m,
\\ p_{n-m+j-1}'(x_i), & m+1\leq i\leq 2m
\end{cases}\right)\\
 = \frac{1}{\Delta(x)^2}\sum_{0\leq k_m<\dots <k_1\leq n-1}\frac{(\det_{1\leq i,j\leq m}(p_{k_i}(x_j)))^2}{\prod_{i=1}^m\|p_{k_i}\|^2},
\end{gather}
\end{subequations}
where the diagonal entries in \eqref{cd} are interpreted as the limit
$$ \lim_{ x_j\rightarrow x_i}\frac{p_n(x_i)p_{n-1}(x_j)-p_{n-1}(x_i)p_n(x_j)}{x_i-x_j}=p_{n}'(x_i)p_{n-1}(x_i)-p_{n-1}'(x_i)p_n(x_i). $$
 To recover  Milne's Schur function expansions we need another identity, namely \cite[Proposition 1.8]{r2},
\begin{equation}\label{cse}C(x)=\sum_{\substack{0\leq \lambda_m\leq\dots\leq \lambda_1\leq n-m\\0\leq \mu_m\leq\dots\leq \mu_1\leq n-m}}\frac{(-1)^{\sum_{i=1}^m(\lambda_i+\mu_i)}}{\prod_{i=1}^n\| p_{i-1}\|^2}\det_{i\in [n]\setminus S,j\in [n]\setminus T}(c_{i+j-2})\,s_\lambda(x)s_\mu(x), \end{equation}
where $c_k=\mu(x^k)$, and where
$$S=\{\lambda_k+m+1-k;\,1\leq k\leq m\},\qquad
T=\{\mu_k+m+1-k;\,1\leq k\leq m\}.$$

An important relation 
 between  Hankel determinants and correlation functions follows 
from  Lemma \ref{hdl} upon replacing
 $\mu$  by $\mu+\lambda$, where we think of $\mu$ as  positive and $\lambda$ as  arbitrary. 
Since the integrand is symmetric, we may write
\begin{multline*}
\det_{1\leq i,j\leq m}\left(\mu(x^{i+j-2})+\lambda(x^{i+j-2})\right)\\
\begin{split}&=\frac 1{m!}\int \Delta(x_1,\dots,x_m)^2\, d(\mu+\lambda)(x_1)\dotsm d(\mu+\lambda)(x_m)
\\
&=\sum_{s=0}^m\frac{1}{s!(m-s)!}\int \Delta(x_1,\dots,x_m)^2\,d\lambda(x_1)\dotsm d\lambda (x_s)d\mu(x_{s+1})\dotsm d\mu(x_m),
\end{split} \end{multline*}
where the integral over the last $m-s$ variables is a correlation function. We thus arrive at the following result, which is a standard tool of random matrix theory; see e.g.\ \cite[\S 2]{j}.

\begin{lemma}\label{hcl}
In the notation above, 
\begin{multline*}
\det_{1\leq i,j\leq m}\left(\mu(x^{i+j-2})+\lambda(x^{i+j-2})\right)\\
=\prod_{i=1}^m\|p_{i-1}\|^2\sum_{s=0}^m \frac 1{s!} \int \Delta(x_1,\dots,x_s)^2\, C_s^m(x_1,\dots,x_s)\,d\lambda(x_1)\dotsm d\lambda(x_s).\end{multline*}
\end{lemma}

\subsection{Milne's sums of squares formulas}\label{sss}

By Lemma \ref{hcl}, the Hankel determinants of Corollary \ref{mhd} may be expressed in terms of correlation functions for continuous dual Hahn polynomials.
If we let  $C_{m}^{n,\ep}$ denote the correlation function defined by choosing  $\mu=\mu_\ep$ in \eqref{kcf}, then 
applying Lemma \ref{hcl}  with $\mu=\mu_0$ and $\lambda=\lambda_0$ yields, using also \eqref{en},
$$\square(q)^{4m^2}=\sum_{s=0}^m\frac{4^s}{s!}\sum_{k_1,\dots,k_s=1}^\infty 
\prod_{i=1}^s \frac{q^{k_i}k_i}{1+(-q)^{k_i}}\prod_{1\leq i<j\leq s}(k_i^2-k_j^2)^2 \,C_s^{m,0}(-k_1^2,\dots,-k_s^2).$$
It is understood that the sum over $k_1,\dots,k_s$ equals $1$ if $s=0$.
Similarly as in Theorem \ref{adf}, we may restrict the summation to $k_1>\dots>k_s$ if we multiply by $s!$, thus obtaining the first half of Corollary \ref{gcc}. The second half follows similarly, choosing  $\mu=\mu_1$ and $\lambda=\lambda_1$
in Lemma \ref{hcl}.

\begin{corollary}
\label{gcc}
One has
$$\square(q)^{4m^2}=\sum_{s=0}^m 4^s\sum_{k_1>\dots>k_s\geq 1}
\prod_{i=1}^s \frac{q^{k_i}k_i}{1+(-q)^{k_i}}\prod_{1\leq i<j\leq s}(k_i^2-k_j^2)^2 \,C_s^{m,0}(-k_1^2,\dots,-k_s^2),$$
\begin{multline*}\square(q)^{4m(m+1)}\\
=\sum_{s=0}^m 4^s\sum_{k_1>\dots>k_s\geq 1}
\prod_{i=1}^s \frac{q^{k_i}k_i^3}{1-(-q)^{k_i}}\prod_{1\leq i<j\leq s}(k_i^2-k_j^2)^2 \,C_s^{m,1}(-k_1^2,\dots,-k_s^2).
\end{multline*}
\end{corollary}

Expanding
$$\frac{q^k}{1+(-q)^k}=\sum_{l=1}^\infty (-1)^{(k-1)(l-1)}q^{kl},\qquad
\frac{q^k}{1-(-q)^k}=\sum_{l=1}^\infty (-1)^{k(l-1)}q^{kl},
 $$
gives the following sums of squares formulas. We have excluded the terms with 
 $s=0$, since they only contribute to the trivial coefficients corresponding to
 $n=0$.

\begin{theorem}\label{sst} For $n>0$, 
\begin{equation*}
\begin{split}\square_{4m^2}(n)
&=\sum_{s=1}^m 4^s\sum_{\substack{k_1l_1+\dots+k_sl_s=n\\k_1>\dots>k_s\geq 1\\l_1,\dots,l_s\geq 1}} \prod_{i=1}^s (-1)^{(k_i-1)(l_i-1)}k_i\prod_{1\leq i<j\leq s}(k_i^2-k_j^2)^2\\
&\quad\times C_s^{m,0}(-k_1^2,\dots,-k_s^2),
\end{split}
\end{equation*}
\begin{equation*}\begin{split}
\square_{4m(m+1)}(n)
&=\sum_{s=1}^m 4^s\sum_{\substack{k_1l_1+\dots+k_sl_s=n\\k_1>\dots>k_s\geq 1\\l_1,\dots,l_s\geq 1}} \prod_{i=1}^s (-1)^{k_i(l_i-1)}k_i^3\prod_{1\leq i<j\leq s}(k_i^2-k_j^2)^2 \\
&\quad\times C_s^{m,1}(-k_1^2,\dots,-k_s^2)
.
\end{split}
\end{equation*}
\end{theorem}

Expressing the correlation functions as in  \eqref{cp} gives explicit versions of these  sums of squares formulas in terms of continuous dual Hahn polynomials. Using instead 
 \eqref{cse}, together with the expressions \eqref{muce} for the moments
 and \eqref{en} for the norms, yields the 
 following identities, which are equivalent to \cite[Theorems~7.1 and~7.2]{mi}.

\begin{corollary}[Milne]
\label{mt}
For $n>0$, 
\begin{multline*}
\square_{4m^2}(n)
=\frac{2^{m(2m-1)}}{\prod_{i=1}^{2m-1}i!}\sum_{s=1}^m 4^s\sum_{\substack{k_1l_1+\dots+k_sl_s=n\\k_1>\dots>k_s\geq 1\\l_1,\dots,l_s\geq 1}} \prod_{i=1}^s (-1)^{(k_i-1)(l_i-1)}k_i\prod_{1\leq i<j\leq s}(k_i^2-k_j^2)^2\\
\times
\sum_{\substack{0\leq \lambda_s\leq\dots\leq \lambda_1\leq m-s\\0\leq \mu_s\leq\dots\leq \mu_1\leq m-s}}\det_{i\in [m]\setminus S,j\in [m]\setminus T}(t_{i+j-1})\,s_\lambda(k_1^2,\dots,k_s^2)s_\mu(k_1^2,\dots,k_s^2),
\end{multline*}
\begin{multline*}
\square_{4m(m+1)}(n)
=\frac{2^{m(2m+1)}}{\prod_{i=1}^{2m}i!}\sum_{s=1}^m 4^s\sum_{\substack{k_1l_1+\dots+k_sl_s=n\\k_1>\dots>k_s\geq 1\\l_1,\dots,l_s\geq 1}} \prod_{i=1}^s (-1)^{k_i(l_i-1)}k_i^3\prod_{1\leq i<j\leq s}(k_i^2-k_j^2)^2\\
\times
\sum_{\substack{0\leq \lambda_s\leq\dots\leq \lambda_1\leq m-s\\0\leq \mu_s\leq\dots\leq \mu_1\leq m-s}}\det_{i\in [m]\setminus S,j\in [m]\setminus T}(t_{i+j})\,s_\lambda(k_1^2,\dots,k_s^2)s_\mu(k_1^2,\dots,k_s^2),
\end{multline*}
where $t_k$ is as in \eqref{cb} and 
$$S=\{\lambda_k+s+1-k;\,1\leq k\leq s\},\qquad
T=\{\mu_k+s+1-k;\,1\leq k\leq s\}.$$
\end{corollary}

\section{Relation to Schur $Q$-polynomials}

In \cite{r3}, we showed that the correlation functions  appearing in Theorem \ref{sst} also arise from  generalized Schur $Q$-polynomials by specializing all variables to $1$. In this section,  we  give a slightly different proof of Milne's formulas, where Schur $Q$-polynomials appear naturally.

\subsection{Schur $Q$-polynomials}

When $m\leq n$ and $\lambda\in\mathbb Z^m$, we write
\begin{equation}\label{pdef}Q_\lambda(x_1,\dots,x_n)=2^m\sum_{\sigma\in
  S_n/S_{n-m}}\sigma\Bigg(x_1^{\lambda_1}\dotsm
  x_m^{\lambda_m}\prod_{\substack{1\leq i\leq m\\1\leq i<j\leq
  n}}\frac{x_i+x_j}{x_i-x_j}\Bigg).\end{equation}
Here, $S_n$ acts by permuting the variables $x_1,\dots,x_n$ and
$S_{n-m}$ is the  subgroup acting on $x_{m+1},\dots,x_n$.
When $\lambda_1\geq\lambda_2\geq\dots\geq\lambda_m> 0$, this is a standard definition of Schur $Q$-polynomials \cite{m}.  Note that there is no analogue of \eqref{st}, so the  case when some $\lambda_i$ are negative  cannot immediately be reduced to the  classical situation.

We will not work with \eqref{pdef} directly, but rather with the alternative expression \cite[Lemma A.1]{r3}
\begin{multline}\label{ne}Q_\lambda(x_1,\dots,x_n)\\
=\frac{2^{m-k}}{ k!}\prod_{1\leq i<j\leq n}\frac{x_i+x_j}{x_i-x_j}
\sum_{\sigma\in S_n}\sgn(\sigma)\prod_{i=1}^m x_{\sigma(i)}^{\lambda_i}\prod_{i=1}^k\frac{x_{\sigma(m+2i-1)}-x_{\sigma(m+2i)}}{x_{\sigma(m+2i-1)}+x_{\sigma(m+2i)}},
\end{multline}
where $k$ is the integral part of $(n-m)/2$. In analogy with \eqref{sf}, the right-hand side can be rewritten as a quotient of two pfaffians \cite[Eq.\ (A12)]{n}; however, the form given here is more convenient for our purposes.

One of the main results of \cite{r3} is the identification of $Q_\lambda(1,q,\dots,q^{n-1})$ with  generalized Christoffel--Darboux kernels for
 continuous $q$-Jacobi polynomials. In the case when $q=1$ and
$\lambda=(\lambda_1,\dots,\lambda_m,-\lambda_m,\dots,-\lambda_1)$,  the continuous $q$-Jacobi polynomials reduce  to continuous dual Hahn polynomials and
the generalized Christoffel--Darboux kernels to  correlation functions. 
To be precise,  noting that the kernels $K_m^{n,\ep}$ featuring in \cite{r3} are related to the kernels 
$C_m^{n,\ep}$ introduced in \S \ref{sss} through
$$C_m^{n,\ep}(x_1,\dots,x_m)=K_m^{n,\ep}(x_1,\dots,x_m,x_1,\dots,x_m), $$
it follows from \cite[Corollary 5.11]{r3} that
\begin{subequations}\label{qc}
\begin{multline}Q_{(\lambda_1,\dots,\lambda_m,-\lambda_m,\dots,-\lambda_1)}(1^{2n})\\
=8^m\prod_{i=1}^{m}\lambda_i\prod_{1\leq i<j\leq m}(\lambda_i^2-\lambda_j^2)^2
\,C_m^{n,0}(-\lambda_1^2,\dots,-\lambda_m^2),
 \end{multline}
\begin{multline}
Q_{(\lambda_1,\dots,\lambda_m,-\lambda_m,\dots,-\lambda_1)}(1^{2n+1})\\
=(-1)^m8^m\prod_{i=1}^{m}\lambda_i^3\prod_{1\leq i<j\leq m}(\lambda_i^2-\lambda_j^2)^2
\,C_m^{n,1}(-\lambda_1^2,\dots,-\lambda_m^2).
 \end{multline}
\end{subequations}

\subsection{A second proof of Milne's formulas}

We  return to the pfaffian evaluations  \eqref{epi}. In \S \ref{fps}, we handled the singularities at $x_i=-x_j$ by introducing Abel means. We will now use a slightly different approach, removing the singularity by subtracting a rational function. This  corresponds precisely to the previous decomposition $\mu_\ep=\lambda_\ep+\nu_\ep$ of functionals, so the two ideas are in fact intimately related.

\begin{lemma}\label{lel}
For $q<|x|<q^{-1}$, 
\begin{subequations}\label{ta}
\begin{equation}\label{tea}\frac{(q)_\infty^2\theta(x)}{(-q)_\infty^2\theta(-x)}=
 \frac{1-x}{1+x}+2\sum_{k=1}^\infty
 \frac{(-q)^k(x^{-k}-x^{k})}{1+q^k},\end{equation}
\begin{equation}\label{toa}1+2x\frac{\theta'(-x)}{\theta(-x)}=\frac{1-x}{1+x}+2\sum_{k=1}^\infty\frac{(-q)^k(x^k-x^{-k})}{1-q^k}. \end{equation}
\end{subequations}
\end{lemma}

\begin{proof}
Subtracting \eqref{rfe} from \eqref{rpa} gives  \eqref{tea} 
for  $q<|x|<1$. By analytic continuation, it holds also for $1\leq |x|<q^{-1}$. 
The equation \eqref{toa} follows similarly  from \eqref{ol}.
\end{proof}

Applying \eqref{ta} to
 the matrix elements in \eqref{epi} leads to pfaffians of the  form $\pfaff(A-A^t+B)$. The following   lemma will be useful.

\begin{lemma}\label{psl}
One has
\begin{multline*}\pfaff_{1\leq i,j\leq m}(a_{ij}-a_{ji}+b_{ij})\\
=\frac{1}{M!}\sum_{s=0}^{M}\frac{1}{2^{{M}-s}}\binom {M}s \sum_{\sigma\in S_m}\sgn(\sigma)\prod_{i=1}^s a_{\sigma(2i-1),\sigma(2i)}\prod_{i=s+1}^{M} b_{\sigma(2i-1),\sigma(2i)},
 \end{multline*}
where $M$ is the integral part of $m/2$  and $(b_{ij})_{i,j=1}^m$ is skew-symmetric.
\end{lemma}

\begin{proof}
Consider first
\begin{eqnarray*}\pfaff_{1\leq i,j\leq m}(a_{ij}+b_{ij})&=&\frac{1}{2^{M}{M}!} \sum_{\tau\in S_m}\sgn(\tau)\prod_{i=1}^{M} (a_{\tau(2i-1),\tau(2i)}+ b_{\tau(2i-1),\tau(2i)})\\
&=&\frac{1}{2^{M}{M}!}\sum_{S\subseteq[{M}]}\,\sum_{\tau\in S_m}\sgn(\tau)\prod_{i\in S} a_{\tau(2i-1),\tau(2i)}\prod_{i\in[{M}]\setminus S} b_{\tau(2i-1),\tau(2i)},
 \end{eqnarray*}
where $[M]=\{1,2,\dots,M\}$.
Choose $\sigma\in S_m$ so that
$$\prod_{i\in S} a_{\tau(2i-1),\tau(2i)}\prod_{i\in[M]\setminus S} b_{\tau(2i-1),\tau(2i)}=\prod_{i=1}^s a_{\sigma(2i-1),\sigma(2i)}\prod_{i=s+1}^M b_{\sigma(2i-1),\sigma(2i)},
$$
where $s=|S|$. This can be done so that, for fixed $S$, $\tau\mapsto \sigma$ is a bijection. Moreover, $\sgn(\sigma)=\sgn(\tau)$.   Rewriting the sum in terms of $\sigma$ and $s$ yields
$$\pfaff_{1\leq i,j\leq m}(a_{ij}+b_{ij})\\
=\frac{1}{2^MM!}\sum_{s=0}^M\binom Ms \sum_{\sigma\in S_m}\sgn(\sigma)\prod_{i=1}^s a_{\sigma(2i-1),\sigma(2i)}\prod_{i=s+1}^M b_{\sigma(2i-1),\sigma(2i)}.
 $$
Replacing $a_{ij}$ by $a_{ij}-a_{ji}$, the sign changes obtained when expanding the product
$\prod_{i=1}^s (a_{\sigma(2i-1),\sigma(2i)}-a_{\sigma(2i),\sigma(2i-1)})$ are compensated by the factor $\sgn(\sigma)$, thus yielding the desired result.
\end{proof}

If $a_{ij}=\sum_{k=1}^\infty c_k(x_i/x_j)^k$ and $b_{ij}=(x_i-x_j)/(x_i+x_j)$, Lemma \ref{psl} gives
\begin{multline*}
\pfaff_{1\leq i,j\leq m}\left(\frac{x_i-x_j}{x_i+x_j}+\sum_{k=1}^\infty c_k\left(\left({x_i}/{x_j}\right)^k-\left({x_j}/{x_i}\right)^k\right)\right)\\
=\frac 1{M!}\sum_{s=0}^M\frac 1{2^{M-s}}\binom Ms
\sum_{k_1,\dots,k_s=1}^\infty\, \sum_{\sigma\in S_m}\prod_{i=1}^s c_{k_i}x_{\sigma(2i-1)}^{k_i}
x_{\sigma(2i)}^{-k_i}\prod_{i=s+1}^M\frac{x_{\sigma(2i-1)}-x_{\sigma(2i)}}{x_{\sigma(2i-1)}+x_{\sigma(2i)}},
\end{multline*}
where, as before,  the sum over $k_1,\dots,k_s$ equals $1$ when $s=0$.
By  \eqref{ne}, 
this may be written 
$$\prod_{1\leq i<j\leq m}\frac{x_i-x_j}{x_i+x_j}\sum_{s=0}^M\frac {1}{4^s s!}\sum_{k_1,\dots,k_s=1}^\infty \prod_{i=1}^s c_{k_i}\,Q_{(k_1,-k_1,\dots,k_s,-k_s)}(x_1,\dots,x_m).$$
Similarly as for Theorem \ref{adf}, we  exploit the symmetry in $k_i$ to write the result as follows.

\begin{proposition}\label{sep}
One has
\begin{multline*}\pfaff_{1\leq i,j\leq m}\left(\frac{x_i-x_j}{x_i+x_j}+\sum_{k=1}^\infty c_k\left(\left({x_i}/{x_j}\right)^k-\left({x_j}/{x_i}\right)^k\right)\right)\\
=\prod_{1\leq i<j\leq m}\frac{x_i-x_j}{x_i+x_j}\sum_{s=0}^M\frac{1}{4^s }\sum_{k_1>k_2>\dots>k_s\geq 1} \,\prod_{i=1}^s c_{k_i}\,Q_{(k_1,\dots,k_s,-k_s,\dots,-k_1)}(x_1,\dots,x_m),
\end{multline*}
where $(c_k)_{k=1}^\infty$ is a scalar sequence such that the series converge absolutely and $M$ is the integral part of $m/2$.
\end{proposition}

\begin{remark}
More generally, the same proof gives
\begin{multline*}\pfaff_{1\leq i,j\leq m}\left(\frac{x_i-x_j}{x_i+x_j}+\sum_{k,l=-\infty}^\infty c_{kl}\left(x_i^kx_j^l-x_j^kx_i^l\right)\right)
=\prod_{1\leq i<j\leq m}\frac{x_i-x_j}{x_i+x_j}\\
\times\sum_{s=0}^M\frac {(-1)^{\frac12{s(s-1)}}}{4^s s!}\sum_{k_1,\dots,k_s,l_1,\dots,l_s=-\infty}^\infty \,\prod_{i=1}^s c_{k_il_i}\,Q_{(k_1,\dots,k_s,l_1,\dots,l_s)}(x_1,\dots,x_m).
\end{multline*}
\end{remark}

We  apply Proposition \ref{sep} when $m$ is even and 
$c_k= 2(-q)^k/(1+ q^k)$, and when $m$ is odd and $c_k=- 2(-q)^k/(1- q^k)$. Then, by Lemma \ref{lel}, the pfaffians are of the form \eqref{epi}. We thus arrive at the following result.  

\begin{theorem}\label{mdt}
If $q<|x_j/x_i|<q^{-1}$ for all $i$ and $j$, then
\begin{multline*}\frac{(q)_\infty^{2m}}{(-q)_\infty^{2m}}\prod_{1\leq i<j\leq 2m}\frac{(qx_j/x_i,qx_i/x_j)_\infty}{(-qx_j/x_i,-qx_i/x_j)_\infty}\\
=\sum_{s=0}^m\frac{1}{2^s }\sum_{k_1>\dots>k_s\geq 1} \prod_{i=1}^s \frac{(-q)^{k_i}}{1+q^{k_i}}\,Q_{(k_1,\dots,k_s,-k_s,\dots,-k_1)}(x_1,\dots,x_{2m}),
 \end{multline*}
\begin{multline*}\frac{(q)_\infty^{2m}}{(-q)_\infty^{2m}}\prod_{1\leq i<j\leq 2m+1}\frac{(qx_j/x_i,qx_i/x_j)_\infty}{(-qx_j/x_i,-qx_i/x_j)_\infty}\\
=\sum_{s=0}^m\frac{(-1)^{s}}{2^s }\sum_{k_1>\dots>k_s\geq 1} \prod_{i=1}^s \frac{(-q)^{k_i}}{1-q^{k_i}}\,Q_{(k_1,\dots,k_s,-k_s,\dots,-k_1)}(x_1,\dots,x_{2m+1}).
 \end{multline*}
\end{theorem}


\begin{remark} 
 It would be interesting to find an algebraic interpretation of
 Theorem \ref{mdt}, beyond the link to denominator formulas for queer affine superalgebras via modular duality. Plausibly, the work of Sergeev \cite{se1,se2}, where Schur $Q$-polynomials arise as characters and spherical functions for queer superalgebras, is relevant for finding such a relation. 
\end{remark}

Specializing as before  $x_i\equiv 1$ and replacing $q$ by $-q$, Theorem \ref{mdt} reduces to the following identities. 

\begin{corollary}\label{qss} One has
$$\square(q)^{4m^2}
=\sum_{s=0}^m\frac{1}{2^s }\sum_{k_1>\dots>k_s\geq 1} \prod_{i=1}^s \frac{q^{k_i}}{1+(-q)^{k_i}}\,Q_{(k_1,\dots,k_s,-k_s,\dots,-k_1)}(1^{2m}),
 $$
$$\square(q)^{4m(m+1)}=\sum_{s=0}^m\frac{(-1)^{s}}{2^s }\sum_{k_1>\dots>k_s\geq 1} \prod_{i=1}^s \frac{q^{k_i}}{1-(-q)^{k_i}}\,Q_{(k_1,\dots,k_s,-k_s,\dots,-k_1)}(1^{2m+1}).$$
\end{corollary}

By \eqref{qc}, this is equivalent to Corollary \ref{gcc}, and thus also to Milne's sums of squares formulas. 

\section{A new formula for $2m^2$  squares}
\label{hss}

As we have seen, Milne's formulas for $4m^2$ and $4m(m+1)$ squares arise from the pfaffian evaluations \eqref{epi} as all variables $x_i\rightarrow 1$. In view of the results for triangular numbers in \cite{r1}, one would expect that more general formulas for $4m^2/d$ squares, when $d\mid 2m$, and $4m(m+1)/d$ squares, when $d\mid 2m$ or $d\mid 2(m+1)$, may be obtained by letting $x_i$ tend to suitable fractional powers of $q$. We have made some preliminary investigations in this direction, but, 
 unless one can simplify the arguments, the  results seem very complicated. For the benefit of the interested reader, we 
state without proof a particularly  simple special case, a $2m^2$ squares formula which we consider to be the natural analogue of \eqref{hti}. It 
 can be obtained from \eqref{eep} in the limit when  $x_1,\dots,x_m\rightarrow 1$, $x_{m+1},\dots,x_{2m}\rightarrow \sqrt q$. 

To state the result we need to introduce the Schur-type 
polynomials 
$$\mathbb P_{n^m}^{(\ep)}(x_1,\dots,x_m)=\frac{\det_{1\leq i,j\leq m}(p_{n+j-1}^{(\ep)}(x_i))}{\Delta(x)},$$
where $p_k^{(\ep)}$ are as in \S \ref{cdhs}. They generalize the correlation functions of \S \ref{sss}, since
 $$C_m^{n,\ep}(x_1,\dots,x_m)=\frac 1{\prod_{i=1}^m\|p_{n-i}^{(\ep)}\|^2}
 \,\mathbb P_{(n-m)^{2m}}^{(\ep)}(x_1,\dots,x_m,x_1,\dots,x_m).$$
Moreover, if
  $n-m=2k+\ep$, with $k$ a
non-negative integer and $\ep\in\{0,1\}$, we have the following generalization of \eqref{qc}
 \cite[Corollary~5.11]{r3}:
 \begin{multline*}Q_{(\lambda_1,\dots,\lambda_m)}(1^n)\\
 =\frac{2^{\frac 12 m(2n+1-m)}(-1)^{km }}{
 \prod_{i=1}^m{(n-i)!}}
 \prod_{i=1}^m\lambda_i^\ep\prod_{1\leq i<j\leq m}(\lambda_i-\lambda_j)
 \,\mathbb P^{(\ep)}_{k^m}(-\lambda_1^2,\dots,-\lambda_m^2).
  \end{multline*}

\begin{theorem}\label{hsf}
One has
\begin{multline*}
\square_{2m^2}(n)=
\sum_{\substack{s_0,s_1,s_2\geq 0\\ s_0+2s_1\leq m,\,s_2\leq s_1}}
(-1)^{(m+1)(s_1+s_2)}(2-\delta_{s_1s_2}) \\
\begin{split}&\times\frac{2^{(s_0+s_1+s_2)(2m+1)-\frac 12((s_0+2s_1)^2+(s_0+2s_2)^2)}}{\prod_{i=1}^{s_0+2s_1}(m-i)!\prod_{i=1}^{s_0+2s_2}(m-i)!}\\
&\times\sum_{\substack{k_1>\dots>k_{t_0}\geq 1,\,k_{t_0+1}>\dots>k_{s_0}\geq 1,\,l_1,\dots,l_{s_0} \text{ \emph{odd positive}}\\
k_1'>\dots>k_{s_1}'\geq 1,\,l_1',\dots,l_{s_1}' \text{ \emph{even positive}}\\
k_1''>\dots>k_{s_2}''\geq 1,\,l_1'',\dots,l_{s_2}'' \text{ \emph{even positive}}\\
k_1l_1+\dots+k_{s_0}l_{s_0}+k_1'l_1'+\dots+k_{s_1}'l_{s_1}'+k_1''l_1''+\dots+k_{s_2}''l_{s_2}''=n}}\,\prod_{i=1}^{s_0}(-1)^{\frac 12(l_i-1)}\prod_{i=1}^{t_0}k_i^{2+2\ep}\prod_{i=t_0+1}^{s_0}k_i^{2\ep}
\\
&\times\prod_{i=1}^{s_1}(-1)^{k_i'+\frac 12 l_i'}(k_i')^{1+2\ep}
\prod_{i=1}^{s_2}(-1)^{k_i''+\frac 12 l_i''}(k_i'')^{1+2\ep}
\\
&\times\prod_{\substack{1\leq i<j\leq t_0\\ \text{\emph{or }} t_0+1\leq i<j\leq s_0}}(k_j^2-k_i^2)^2\prod_{1\leq i\leq s_0,\,1\leq j\leq s_1}
(k_i^2-(k_j')^2)\prod_{1\leq i\leq s_0,\,1\leq j\leq s_2}
(k_i^2-(k_j'')^2)\\
&\times \prod_{1\leq i<j\leq s_1}((k_j')^2-(k_i')^2)^2\prod_{1\leq i<j\leq s_2}((k_j'')^2-(k_i'')^2)^2\\
&\times \mathbb P^{(\ep)}_{t_1^{s_0+2s_1}}(-k_1^2,\dots,-k_{s_0}^2,-(k_1')^2,-(k_1')^2,\dots,-(k_{s_1}')^2,-(k_{s_1}')^2)\\
&\times \mathbb P^{(\ep)}_{t_2^{s_0+2s_2}}(-k_1^2,\dots,-k_{s_0}^2,-(k_1'')^2,-(k_1'')^2,\dots,-(k_{s_2}'')^2,-(k_{s_2}'')^2)
,
\end{split}\end{multline*}
where  $\ep=0$ if $m-s_0$ is even and $1$ else, and where
$t_0$, $t_1$ and $t_2$ denote the integral part of $s_0/2$, $(m-s_0-2s_1)/2$ and $(m-s_0-2s_2)/2$, respectively.
\end{theorem}

As examples, we work out the cases $m=1,\,2,\, 3$. 
We observe that the term $s=(s_0,s_1,s_2)=(0,0,0)$ only contributes to the trivial case when $n=0$, so we assume that $n>0$. Note also that when $m=1$ or $2$, all Schur-type polynomials $\mathbb P_{k^l}$ that appear have either  $k=0$ or $l=0$, and thus equal $1$.

 \emph{2 squares}: When $m=1$, the outer sum  has only one non-trivial term, $s=(1,0,0)$. We recover the two squares formula \eqref{s2}. 

 \emph{8 squares}: When $m=2$, we obtain an eight squares formula, which is more complicated than \eqref{s8} but still interesting to discuss. We  have four non-trivial terms,  $s=(0,1,0)$, $(1,0,0)$, $(0,1,1)$ and $(2,0,0)$, giving
\begin{multline*}
\square_8(n)=16\sum_{k_1'l_1'=n,\,l_1' \text{ even}}(-1)^{k_1'-1+\frac 12 l_1'}k_1'
+ 16\sum_{k_1l_1=n,\,l_1 \text{ odd}}(-1)^{\frac 12 (l_1-1)}k_1^2\\
+64 \sum_{\substack{k_1'l_1'+k_1''l_1''=n\\ l_1' \text{ and } l_1'' \text{ even}}}(-1)^{k_1'+\frac 12l_1'+k_1''+\frac 12l_1''}k_1'k_1''
+64 \sum_{\substack{k_1l_1+k_2l_2=n\\ l_1 \text{ and } l_2 \text{ odd}}}(-1)^{\frac 12(l_1+l_2-2)}k_1^2.
 \end{multline*}
Equivalently, 
\begin{multline*}
\square(q)^8=1+16\sum_{k=1}^\infty\frac{(-1)^kkq^{2k}}{1+q^{2k}}+16\sum_{k=1}^\infty\frac{k^2q^{k}}{1+q^{2k}}\\
+64\left(\sum_{k=1}^\infty\frac{(-1)^kkq^{2k}}{1+q^{2k}}\right)^2
+64 \sum_{k=1}^\infty\frac{k^2q^{k}}{1+q^{2k}}\sum_{k=1}^\infty\frac{q^{k}}{1+q^{2k}}
.\end{multline*}
Computing
the Lambert series  using  \eqref{l2}, \eqref{l4} and \eqref{jl} gives after simplification
$$\square(q)^8=\square(-q^2)^8+16q\triangle(q^2)^4\square(q)^4. $$
Applying
$\square(-q^2)^2=\square(q)\square(-q)$,
which is easily verified either from \eqref{sp} or from the definition,
we recover Jacobi's quartic identity \eqref{jq}.

\emph{18 squares:}  When $m=2$, we have seven non-trivial terms, $s=(0,1,0)$, $(1,0,0)$, $(0,1,0)$, $(1,1,0)$, $(2,0,0)$, $(1,1,1)$ and $(3,0,0)$. 
The second and fourth term involve the polynomial $\mathbb P_{1^1}^{(0)}(x)=p_1^{(0)}(x)=x-\frac 12$, while all other Schur-type factors are trivial. 
We thus obtain
\begin{multline*}\begin{split}
\square_{18}(n)&=2^5\sum_{k_1'l_1'=n,\,l_1' \text{ even}}(-1)^{k_1'+\frac 12 l_1'} (k_1')^3\\
&\quad
+ 2^4\sum_{k_1l_1=n,\,l_1 \text{ odd}}(-1)^{\frac 12 (l_1-1)}\left(k_1^2+\frac 12\right)^2\\
&\quad
+2^8 \sum_{\substack{k_1'l_1'+k_1''l_1''=n\\ l_1' \text{ and } l_1'' \text{ even}}}(-1)^{k_1'+\frac 12l_1'+k_1''+\frac 12l_1''}(k_1')^3(k_1'')^3\\
&\quad+2^8 \sum_{\substack{k_1l_1+k_1'l_1'=n\\ l_1 \text{ odd and } l_1' \text{ even}}}(-1)^{\frac 12(l_1-1)+k_1'-1+\frac 12l_1'}\left(k_1^2+\frac 12\right)k_1'\\
&\quad+2^8 \sum_{\substack{k_1l_1+k_2l_2=n\\ l_1 \text{ and } l_2 \text{ odd}}}(-1)^{\frac 12(l_1+l_2-2)}k_1^4k_2^2\\
&\quad+2^{10} \sum_{\substack{k_1l_1+k_1'l_1'+k_1''l_1''=n\\ l_1 \text{ odd},\, l_1', l_1'' \text{ even}}}(-1)^{\frac 12(l_1-1)+k_1'+\frac 12l_1'+k_1''+\frac 12l_1''}k_1'k_1''\\
&\quad+2^{10} \sum_{\substack{k_1l_1+k_2l_2+k_3l_3=n\\ l_1,l_2,l_3 \text{ odd}}}(-1)^{\frac 12(l_1+l_2+l_3-3)}k_1^2(k_2^2-k_3^2)^2.
\end{split} \end{multline*}

\section*{Appendix. Ono's sums of squares formulas}
\setcounter{section}{1}
\renewcommand{\thesection}{\Alph{section}}
\setcounter{equation}{0}

In this Appendix we show that Ono's sums of squares formulas \cite{o} are equivalent to those of  Milne. Using the notation of Ono, we define
 $A_m^\pm$  by
$$\prod_{i=1}^m x_i\prod_{1\leq i<j\leq m}(x_j^2-x_i^2)^2=\sum_{\lambda=(a_1,\dots,a_m)}A_m^+(\lambda)\,x_1^{a_1}\dotsm x_m^{a_m}, $$
$$\prod_{i=1}^m x_i^3\prod_{1\leq i<j\leq m}(x_j^2-x_i^2)^2=\sum_{\lambda=(a_1,\dots,a_m)}A_m^-(\lambda)\,x_1^{a_1}\dotsm x_m^{a_m}. $$
Moreover, we need the modular forms
\begin{multline*}E^+(2k)= 2^{4k-1}\left(\frac{(-1)^k|B_{2k}|}{4k}+\sum_{n=1}^\infty\sigma_{2k-1}(n)q^{4n} \right)\\
-2^{2k-1}\left(\frac{(-1)^k|B_{2k}|}{4k}+\sum_{n=1}^\infty\sigma_{2k-1}(n)q^{n} \right), \end{multline*} 
\begin{multline*}E^-(2k)= 2^{2k}\left(\frac{(-1)^k|B_{2k}|}{4k}+\sum_{n=1}^\infty\sigma_{2k-1}(n)q^{2n} \right)\\
-\left(\frac{(-1)^k|B_{2k}|}{4k}+\sum_{n=1}^\infty(-1)^n\sigma_{2k-1}(n)q^{n} \right), 
\end{multline*} 
where $\sigma_k(n)=\sum_{d\mid n}d^k$. Compared to \cite{o}, we have replaced $q$ by $-q$, and used that $B_{2i}=(-1)^{i-1}|B_{2i}|$. 
In generating function form, Ono's result then takes the following form \cite[Theorem 1]{o}.

\begin{theorem}[Ono] \label{ot} In the notation above,
$$\square(q)^{4m^2}=\frac{(-1)^m 4^m}{m!\prod_{i=1}^{2m-1}i! }
\sum_{\lambda=(a_1,\dots,a_m)}A_m^+(\lambda) E^+(a_1+1)\dotsm E^+(a_m+1),
 $$
$$\square(q)^{4m(m+1)}=\frac{ 2^{2m^2+3m}}{m!\prod_{i=1}^{2m}i! }
\sum_{\lambda=(a_1,\dots,a_m)}A_m^-(\lambda) E^-(a_1+1)\dotsm E^-(a_m+1).
 $$
\end{theorem}

We claim that Theorem \ref{ot} is equivalent to  Corollary \ref{mhd}. To see this we observe that,  in the notation of Lemma \ref{hdl},
$$A_m^+(a_1,\dots,a_m)=C(k_1,\dots,k_m),\qquad a_i=2k_i+ 1, $$
$$A_m^-(a_1,\dots,a_m)=C(k_1,\dots,k_m),\qquad a_i=2k_i+ 3. $$
Thus, Theorem \ref{ot} can be written in Hankel determinant form as
$$\square(q)^{4m^2}=\frac{(-1)^m 4^m}{\prod_{i=1}^{2m-1}i! }
\det_{1\leq i,j\leq m}(E^+(2i+2j-2)),$$
$$\square(q)^{4m(m+1)}=\frac{ 2^{2m^2+3m}}{\prod_{i=1}^{2m}i! }
\det_{1\leq i,j\leq m}(E^-(2i+2j)).
 $$
 Our claim would now would follow from the identities
$$E^+(2k)=(-1)^k 2^{2k-3}\left(\frac{(4^k-1)|B_{2k}|}{k}+4(-1)^{k+1}\sum_{n=1}^\infty\frac{q^nn^{2k-1}}{1+(-q)^n}\right), $$
$$E^+(2k)=\frac{(-1)^k}4 \left(\frac{(4^k-1)|B_{2k}|}{k}+4(-1)^{k}\sum_{n=1}^\infty\frac{q^nn^{2k-1}}{1-(-q)^n}\right). $$
Cancelling
the terms involving Bernoulli numbers, we are reduced to proving that
$$\sum_{n=1}^\infty\sigma_{2k-1}(n)q^{n}-4^k\sum_{n=1}^\infty \sigma_{2k-1}(n)q^{4n}= \sum_{n=1}^\infty\frac{q^nn^{2k-1}}{1+(-q)^n},$$
$$4^k\sum_{n=1}^\infty \sigma_{2k-1}(n)q^{2n}+\sum_{n=1}^\infty(-1)^{n-1}\sigma_{2k-1}(n)q^{n}= \sum_{n=1}^\infty\frac{q^nn^{2k-1}}{1-(-q)^n},$$
or, equivalently, replacing $q$ by $-q$ in   the second identity, to the elementary identities
\begin{equation}\label{oe}\sum_{l,m\geq 1}l^{2k-1}q^{lm}-2\sum_{l,m\geq 1}(2l)^{2k-1}q^{4lm}
=\sum_{l,m\geq 1}(-1)^{(l-1)(m-1)}l^{2k-1}q^{lm},
 \end{equation}
$$2\sum_{l,m\geq 1}(2l)^{2k-1}q^{2lm}-\sum_{l,m\geq 1}l^{2k-1}q^{lm}
=\sum_{l,m\geq 1}(-1)^{l}l^{2k-1}q^{lm}.
 $$


\begin{thebibliography}{99}
\bibitem[CK]{ck} H.\ H.\ Chan and C.\ Krattenthaler, \emph{Recent progress in the study of representations of integers as sums of squares}, Bull.\ London Math.\ Soc.\ 37 (2005), 818--826.
\bibitem[F]{f} F.\ G.\ Frobenius, \emph{\"Uber die elliptischen
    Functionen zweiter Art},  J.\ Reine Angew.\ Math.\ 93 (1882), 53--68.
\bibitem[FS]{fs}  F.\ G.\ Frobenius and L.\ Stickelberger, \emph{Zur
    Theorie der elliptischen Functionen}, J.\ Reine Angew.\ Math.\ 83 (1877), 175--179.
\bibitem[GR]{gr} G.\ Gasper and M.\ Rahman,  Basic Hypergeometric Series, $2^{\text{nd}}$ ed.,
Cambridge University Press, Cambridge, 2004.
\bibitem[GM]{gm} J.\ Getz and K.\ Mahlburg,
\emph{Partition identities and a theorem of Zagier},
J.\ Combin.\ Theory Ser.\ A 100 (2002),  27--43.
\bibitem[I]{i} M.\ E.\ H.\ Ismail, Classical and Quantum Orthogonal Polynomials in One Variable, Cambridge University Press, Cambridge, 2005.
\bibitem[J]{j} K.\ Johansson, \emph{Random matrices and determinantal processes}, math-ph/0510038.
\bibitem[KW]{kw} V.\ G.\ Kac and M.\ Wakimoto, \emph{Integrable
    highest weight modules over affine superalgebras and number
    theory}, in  Lie Theory and Geometry, pp.\
415--456, Progr.\ Math.\
    123, Birkh\"auser, Boston, 1994.
\bibitem[KS]{ks}  R.\ Koekoek and R.\ F.\ Swarttouw,  The Askey-Scheme of
Hypergeometric Orthogonal Polynomials and its $q$-Analogue,
Delft University of Technology, 1998.
\bibitem[K]{ko} W.\ K\"onig, \emph{Orthogonal polynomial ensembles in probability theory},  Probab.\ Surv.\  2  (2005), 385--447.
\bibitem[LY]{ly} L.\ Long and Y.\ Yang, \emph{A short proof of Milne's formulas for sums of integer squares}, Int.\ J.\ Number Theory~1  (2005),   533--551.
\bibitem[Ma]{m} I.\ G.\ Macdonald, Symmetric Functions and Hall
  Polynomials, $2^{\text{nd}}$ ed., Oxford University Press, Oxford,
  1995.
\bibitem[M1]{mi1} S.\ C.\ Milne,
\emph{New infinite families of exact sums of squares formulas, Jacobi
  elliptic  functions,  and Ramanujan's tau function},
Proc.\ Nat.\ Acad.\ Sci.\ U.S.A. 93 (1996), 15004--15008.
\bibitem[M2]{mi2} S.\ C.\ Milne,
\emph{New infinite families of exact sums of squares formulas,
Jacobi elliptic functions, and Ramanujan's tau function}, 
in  Formal Power Series and Algebraic Combinatorics, 9th Conference, Volume 3,
pp.\ 403--417, Universit\"at Wien, 1997.
\bibitem[M3]{mi}  S.\ C.\ Milne,
\emph{Infinite families of exact sums of squares formulas, Jacobi elliptic functions, continued fractions, and Schur functions},
Ramanujan J.\ 6 (2002),  7--149.
\bibitem[N]{n} J.\ J.\ C.\ Nimmo, \emph{Hall--Littlewood symmetric functions and the BKP equation}, J.\ Phys.\ A 23 (1990), 751--760.
\bibitem[O]{ok} S.\ Okada, \emph{An elliptic generalization of Schur's
    Pfaffian identity}, Adv.\ Math.\  204 (2006),  530--538. 
\bibitem[On]{o} K.\ Ono, \emph{Representations of integers as sums of squares},
J.\ Number Theory 95 (2002), 253--258.
\bibitem[R1]{r1} H.\ Rosengren, \emph{Sums of triangular numbers from the Frobenius determinant}, Adv.\ Math., to appear.
\bibitem[R2]{r3}  H.\ Rosengren, \emph{Schur $Q$-polynomials, multiple hypergeometric series and enumeration of marked shifted tableaux}, math.CO/0603086.
\bibitem[R3]{r2} H.\ Rosengren, \emph{Multivariable Christoffel--Darboux kernels and characteristic polynomials of random hermitian matrices}, math.CA/0606391.
\bibitem[Sc]{s} I.\ Schur,  \emph{\"Uber die Darstellung der
    symmetrischen und der alternierenden Gruppe durch gebrochene
    lineare Substitutionen},  J.\ Reine Angew.\ Math.\ 139 (1911),
    155--250.
19--427. 
\bibitem[S1]{se1} A.\ Sergeev, \emph{Tensor algebra of the identity representation as a module over the Lie superalgebras $GL(n,m)$ and ${\mathcal Q}(n)$,}  Math.\ USSR Sbornik 51 (1985), 419--427. 
\bibitem[S2]{se2} A.\ Sergeev, 
\emph{Projective Schur functions as bispherical functions on certain homogeneous superspaces}, in 
The orbit method in geometry and physics, pp.\ 421--443,
Progr.\ Math.\ 213,
Birkh\"auser, Boston, 2003. 
\bibitem[T]{t} 
E.\ C.\ Titchmarsh, 
The Theory of the Riemann Zeta-Function,  Clarendon Press, Oxford, 1951.  
\bibitem[Z]{z}  D.\ Zagier, \emph{A proof of the Kac--Wakimoto affine
    denominator formula for the strange
series},  Math.\ Res.\ Lett.\  7  (2000),   597--604.
\vskip 3mm
\end{thebibliography}
\end{document}